\documentclass[graybox]{svmult}

\usepackage{type1cm}        
\usepackage{makeidx}         
\usepackage{graphicx}        
\usepackage{multicol}        
\usepackage[bottom]{footmisc}

\usepackage{newtxtext}       %
\usepackage[varvw]{newtxmath}       

\usepackage[numbers]{natbib}
\makeindex             

\newcommand{\Curl}{\ensuremath{\nabla\times}}
\newcommand{\Div}{\ensuremath{\nabla\cdot}}

\newcommand{\Afield}{\ensuremath{\boldsymbol{A}}}
\newcommand{\Bfield}{\ensuremath{\boldsymbol{B}}}

\newcommand{\Efield}{\ensuremath{\boldsymbol{E}}}
\newcommand{\Jfield}{\ensuremath{\boldsymbol{J}}}

\newcommand{\nfield}{\ensuremath{\boldsymbol{n}}}

\newcommand{\zfield}{\ensuremath{\boldsymbol{0}}}

\newcommand{\conductivity}{\ensuremath{\sigma}}
\newcommand{\reluctivity}{\ensuremath{\nu}}
\newcommand{\dualityPairing}[3]{\ensuremath{\big\langle#1,#2\big\rangle_{#3}}}
\newcommand{\dV}{\ensuremath{\operatorname{d}\!\mathrm{V}}}

\newcommand{\ddt}{\ensuremath{\partial t}}

\begin{document}

    \title*{Tree-Cotree-Based IETI-DP for Eddy Current Problems in Time-Domain}
    \author{Mario Mally\orcidID{0009-0000-2685-3392}, Rafael Vázquez\orcidID{0000-0003-1305-6970} and Sebastian Schöps\orcidID{0000-0001-9150-0219}}
    \institute{Mario Mally \at Computational Electromagnetics Group, Technische Universität Darmstadt, 64289 Darmstadt, Germany and Department of Applied Mathematics, Universidade de Santiago de Compostela, 15782 Santiago de Compostela, Spain. \email{mario.mally@tu-darmstadt.de}
    \and Rafael Vázquez \at Department of Applied Mathematics, Universidade de Santiago de Compostela and Galician Centre for Mathematical Research and Technology (CITMAga), 15782 Santiago de Compostela, Spain. \email{rafael.vazquez@usc.es}
    \and Sebastian Schöps \at Computational Electromagnetics Group, Technische Universität Darmstadt, 64289 Darmstadt, Germany \email{sebastian.schoeps@tu-darmstadt.de}
    }
    %
    %
    \maketitle

    \abstract*{For low-frequency electromagnetic problems, where wave-propagation effects can be neglected, eddy current formulations are commonly used as a simplification of the full Maxwell’s equations. In this setup, time-domain simulations, needed to capture transient startup responses or nonlinear behavior, are often computationally expensive. We propose a novel tearing and interconnecting approach for eddy currents in time-domain and investigate its scalability.}

    \abstract{For low-frequency electromagnetic problems, where wave-propagation effects can be neglected, eddy current formulations are commonly used as a simplification of the full Maxwell’s equations. In this setup, time-domain simulations, needed to capture transient startup responses or nonlinear behavior, are often computationally expensive. We propose a novel tearing and interconnecting approach for eddy currents in time-domain and investigate its scalability.}

    \section{Introduction}
    Eddy current formulations are widely employed for low-frequency applications in which the electric energy plays no major role but Ohmic losses and the magnetic energy need to be considered \cite[Sec.~1.2]{Alonso-Rodriguez_2010aa}. Therein, induction heating and levitation/braking devices are presented, but other applications such as induction machines or power loss in transformers \cite{Salon_2023aa} are possible as well. For transient startup responses and nonlinear or even hysteretic behavior, simulations in the time-domain are necessary. These are inherently time-consuming, which, in particular, motivates the use of domain decomposition approaches. Eddy current problems can be approached using different formulations and also various discretization schemes as well as domain decomposition methods. Here, we focus on the $A^{*}$-formulation \cite{Emson_1988aa}, discretize the corresponding weak formulation with isogeometric analysis (IGA) \cite{Cottrell_2009aa} and introduce a dual-primal tearing and interconnecting (IETI-DP) approach following the original FETI-DP idea from \cite{Farhat_2001aa}. Conceptually, we consider our IETI-DP method as an extension of the classical FETI-DP \cite{Farhat_2001aa} to IGA, which allows treating application-oriented problems without introducing geometric approximation errors. In this context, a tree-cotree decomposition is employed to consistently gauge insulating regions, as in \cite{Mally_2025ab} for magnetostatics. A similar FETI-DP investigation was carried out in \cite{Yao_2012ab}, where they use a mixed $A$-$\Phi$-formulation in frequency-domain. In comparison, our work features an extension to the time-domain, investigates scalability and discusses the necessary consistency conditions on interfaces between conducting and insulating regions. This paper is structured as follows. First, in Sec.~\ref{sec:eddy}, the eddy current formulation, including mortaring-based TI coupling, is discussed in the continuous and discrete sense. Then, the dual-primal concept and our tree-cotree approach are explained in Sec.~\ref{sec:DD+TCG}. At last, a numerical experiment for verification and a closer investigation of scalability indicators is carried out in Sec.~\ref{sec:numers}.
    
    \section{Eddy Current Formulations}\label{sec:eddy}
    Following \cite{Alonso-Rodriguez_2010aa}, we assume that the eddy current problem in an open, bounded and simply-connected domain $\Omega\subset\mathbb{R}^3$ is given as
    \begin{align}
        \Curl\left(\reluctivity\Bfield\right) &= \conductivity\Efield + \Jfield, \label{eq:ampere}\\
        \Curl\Efield &= -\ddt\Bfield, \label{eq:faraday}\\
        \Div\Bfield &= 0, \label{eq:gaussM}
    \end{align}
    for every time $t\in\mathcal{I}=(0,T)$. In \eqref{eq:ampere}-\eqref{eq:gaussM}, $\Bfield\colon\Omega\times\mathcal{I}\rightarrow\mathbb{R}^3$ is the magnetic flux density, $\Efield\colon\Omega\times\mathcal{I}\rightarrow\mathbb{R}^3$ the electric field strength and $\Jfield\colon\Omega\times\mathcal{I}\rightarrow\mathbb{R}^3$ a prescribed current density acting as a source. To preserve conciseness, we focus on linear material, but note that nonlinear behavior can be modeled as a simple extension of the approaches discussed in this treatise. Accordingly, the reluctivity $\reluctivity\colon\Omega\rightarrow\mathbb{R}^{+}$ is uniformly bounded, i.e., $0<\reluctivity_{\min}\leq\reluctivity\leq\reluctivity_{\max}$, while the conductivity is only given as $\conductivity\colon\Omega_{\mathrm{C}}\rightarrow\mathbb{R}^{+}$ in open and connected $\Omega_{\mathrm{C}}\subset\Omega$ with $0<\conductivity_{\min}\leq\conductivity\leq\conductivity_{\max}$. In the remaining, insulating part $\Omega_{\mathrm{I}}=\Omega\setminus\overline{\Omega_{\mathrm{C}}}$, we assume $\conductivity=0$ and write $\Gamma=\partial\Omega_{\mathrm{C}}\cap\partial\Omega_{\mathrm{I}}$ for the interfaces shared by conductor and insulator. Let its normal vector $\nfield_{\Gamma}$ point in the direction of $\Omega_{\mathrm{I}}$. Furthermore, we assume that perfect electric conductor boundary conditions (PEC), i.e., $\Efield\times\nfield=\zfield$, on all of $\partial\Omega$, with $\nfield$ being the outward-pointing normal, are given. To solve \eqref{eq:ampere}-\eqref{eq:gaussM}, we use the $A^{*}$-formulation \cite{Emson_1988aa}, for which one defines a vector potential $\Afield\colon\Omega\times\mathcal{I}\rightarrow\mathbb{R}^3$ such that
    \begin{equation}
        \Bfield = \Curl\Afield\qquad\text{and}\qquad\Efield = -\ddt\Afield.\label{eq:vecPot}
    \end{equation}
    Consequently, we obtain the potential equation
    \begin{equation}
        \Curl\left(\reluctivity\Curl\Afield\right) + \conductivity\ddt\Afield = \Jfield,\label{eq:formA*}
    \end{equation}
    by inserting \eqref{eq:vecPot} in \eqref{eq:ampere}. Note that the PEC condition is typically expressed as $\Afield\times\nfield=\zfield$ on $\partial\Omega$ and for all $t\in\mathcal{I}$. This corresponds to a homogeneous Dirichlet boundary condition. Inhomogeneous boundary conditions can be incorporated by using a classical homogenization technique, which we do not discuss in detail here. At last, we assume that an appropriate initial condition $\Afield(\vec{x},t=0) = \Afield_0(\vec{x})$ is given. But this is not enough for \eqref{eq:formA*} to be uniquely solvable. As $\conductivity$ vanishes in $\Omega_{\mathrm{I}}$, a solution can be modified with a gradient field in $\Omega_{\mathrm{I}}$ while still satisfying \eqref{eq:formA*}. We deal with this using a tree-cotree decomposition and dual-primal TI, which is explained in Sec.~\ref{sec:DD+TCG}.

    \subsection{Weak Problem}
    In the following, we use the notation
    \begin{equation*}
        \dualityPairing{\boldsymbol{U}}{\boldsymbol{V}}{\Omega_{i}}=\int_{\Omega_{i}}\boldsymbol{U}\cdot\boldsymbol{V}\dV,~~~\dualityPairing{u}{v}{\Omega_{i}}=\int_{\Omega_{i}}uv\dV,~~~\dualityPairing{\boldsymbol{\eta}}{\boldsymbol{\zeta}}{\Gamma}=\int_{\Gamma}\boldsymbol{\eta}\cdot\boldsymbol{\zeta}\dV,
    \end{equation*}
    for $i\in\{\mathrm{C},\mathrm{I}\}$. Let the classical function spaces $H^1(\Omega_{i})$ and $H(\mathrm{curl},\Omega_{i})$ be given. Then, we can further define the local spaces equipped with Dirichlet boundary conditions as
    \begin{equation*}
        \mathbb{W}_i = \left\{\boldsymbol{U}\in H(\mathrm{curl},\Omega_{i})~~\vert~~\boldsymbol{U}\times\nfield=\zfield~~\text{on}~~\partial\Omega_{i}\cap\partial\Omega\right\}.
    \end{equation*}
    Using these, we can express our mixed mortar formulation for eddy current problems as: Find $\Afield_{\mathrm{C}}\in\mathbb{W}_{\mathrm{C}}$, $\Afield_{\mathrm{I}}\in\mathbb{W}_{\mathrm{I}}$ and $\boldsymbol{\lambda}\in\mathbb{M}$ such that
    \begin{alignat}{2}
        \dualityPairing{\conductivity\ddt\Afield_{\mathrm{C}}}{\boldsymbol{V}_{\mathrm{C}}}{\Omega_{\mathrm{C}}} + a_{\mathrm{C}}\left(\Afield_{\mathrm{C}},\boldsymbol{V}_{\mathrm{C}}\right) + b\left(\boldsymbol{\lambda},\boldsymbol{V}_{\mathrm{C}}\right) &= \dualityPairing{\Jfield_{\mathrm{C}}}{\boldsymbol{V}_{\mathrm{C}}}{\Omega_{\mathrm{C}}},~~&&\forall\boldsymbol{V}_{\mathrm{C}}\in\mathbb{W}_{\mathrm{C}}, \label{eq:saddle1} \\
        a_{\mathrm{I}}\left(\Afield_{\mathrm{I}},\boldsymbol{V}_{\mathrm{I}}\right) - b\left(\boldsymbol{\lambda},\boldsymbol{V}_{\mathrm{I}}\right) &= \dualityPairing{\Jfield_{\mathrm{I}}}{\boldsymbol{V}_{\mathrm{I}}}{\Omega_{\mathrm{I}}},~~&&\forall\boldsymbol{V}_{\mathrm{I}}\in\mathbb{W}_{\mathrm{I}}, \label{eq:saddle2} \\
        b\left(\boldsymbol{\mu},\Afield_{\mathrm{C}} - \Afield_{\mathrm{I}}\right) &= 0,~~&&\forall\boldsymbol{\mu}\in\mathbb{M}, \label{eq:saddle3}
    \end{alignat}
    for every time step $t\in\mathcal{I}$. Note that 
    \begin{equation*}
        a_i\left(\Afield_i,\boldsymbol{V}_i\right) = \dualityPairing{\reluctivity_i\Curl\Afield_i}{\Curl\boldsymbol{V}_i}{\Omega_i}\quad\text{and}\quad
        b\left(\boldsymbol{\lambda},\boldsymbol{V}_i\right) = \dualityPairing{\boldsymbol{\lambda}}{\nfield_{\Gamma}\times\boldsymbol{V}_i\times\nfield_{\Gamma}}{\Gamma}
    \end{equation*}
    for $i\in\{\mathrm{C},\mathrm{I}\}$ and that we followed \cite{Buffa_2020aa} for the underlying structure of \eqref{eq:saddle1}-\eqref{eq:saddle3} and refer to it for more information on the multiplier space $\mathbb{M}$.
    \begin{figure}
        \centering
        \includegraphics[height=2.6cm]{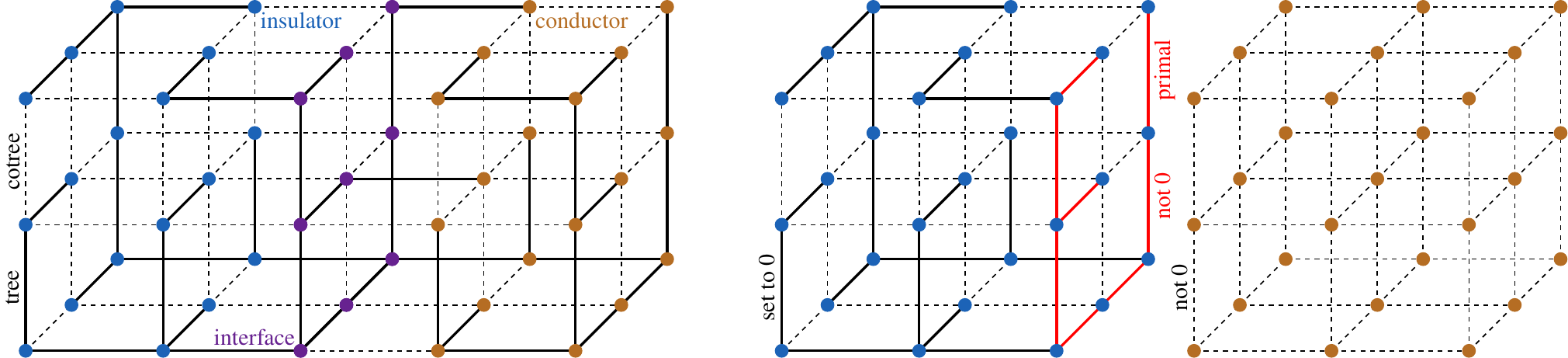}
        \caption{Global tree-cotree splitting (left) and consistent, local gauging (right).}
        \label{fig:tree-cotree}
    \end{figure}

    \subsection{Discrete Problem}
    In the following, we employ a so-called method of lines approach, in which a semi-discrete problem is derived by discretizing in space first before employing a time-stepping scheme. For the discretization, we employ the IGA framework explained in \cite{Vazquez_2016aa}, i.e., we use spline-based, edge-element basis functions $\boldsymbol{w}^{(i)}_j$, for $j\in\{1,\ldots,n_{i}\}$, which span the space $\mathbb{S}_{i}\subset \mathbb{W}_i$. We further require that $\mathbb{S}_{\mathrm{C}}$ and $\mathbb{S}_{\mathrm{I}}$ are conforming, i.e., have pairwise matching basis functions, on $\Gamma$. To approximate $\boldsymbol{\lambda}\in\mathbb{M}$, we assume that a biorthogonal basis
    \begin{equation}
        \operatorname{span}\left(\boldsymbol{\psi}_{j}\right)_{j=1}^{m}\subset \mathbb{M}\quad\text{s.t.}\quad b\left(\boldsymbol{\psi}_{j},\boldsymbol{w}^{(i)}_k\right)=\delta_{jk},~~k\in\mathcal{J}^{(i)}_{\Gamma},j\in\{1,\ldots,m\}\label{eq:biorth}
    \end{equation}
    is given, where $\operatorname{card}(\mathcal{J}^{(i)}_{\Gamma})=m$ for the indices $\mathcal{J}^{(i)}_{\Gamma}$ of all basis functions with non-vanishing support on $\Gamma$. Note that the biorthogonality in \eqref{eq:biorth} is valid for both $\mathbb{S}_{\mathrm{C}}$ and $\mathbb{S}_{\mathrm{I}}$ at the same time, because the spaces are conforming. As a consequence, the coupling reduces to a combination of restriction matrices $\mathbf{R}_{i}$ of size $m\times n_{i}$. These are Boolean matrices and contain only one entry per row which effectively selects certain DOFs. We define the matrices
    \begin{equation*}
        \left(\mathbf{K}_{i}\right)_{jk} = a_i\left(\boldsymbol{w}^{(i)}_k,\boldsymbol{w}^{(i)}_j\right),~~\left(\mathbf{M}_{\mathrm{C}}\right)_{jk} = \dualityPairing{\conductivity\boldsymbol{w}^{(\mathrm{C})}_k}{\boldsymbol{w}^{(\mathrm{C})}_j}{\Omega_{\mathrm{C}}},~~
        \left(\mathbf{j}_{i}\right)_{j} = \dualityPairing{\Jfield_{i}}{\boldsymbol{w}^{(i)}_j}{\Omega_i},
    \end{equation*}
    using $j,k\in\{1,\ldots,n_{i}\}$ and obtain all tools to represent \eqref{eq:saddle1}-\eqref{eq:saddle3} as
    \begin{equation}
        \begin{bmatrix}
            \mathbf{M} & \mathbf{0} \\
            \mathbf{0} & \mathbf{0} \\
        \end{bmatrix}\begin{bmatrix}
            \dot{\mathbf{a}} \\
            \dot{\mathbf{m}} \\
        \end{bmatrix} + \begin{bmatrix}
            \mathbf{K} & \mathbf{B}^{\top} \\
            \mathbf{B} & \mathbf{0} \\
        \end{bmatrix}\begin{bmatrix}
            \mathbf{a} \\
            \mathbf{m} \\
        \end{bmatrix} = \begin{bmatrix}
            \mathbf{j} \\
            \mathbf{0} \\
        \end{bmatrix}.\label{eq:semiDiscDAE}
    \end{equation}
    Therein, we assume $\mathbf{a}(t)\in\mathbb{R}^{n_{\mathrm{C}} + n_{\mathrm{I}}}$ and $\mathbf{m}(t)\in\mathbb{R}^{m}$ for every $t\in\mathcal{I}$ as well as
    \begin{equation*}
        \mathbf{M} = \begin{bmatrix}
            \mathbf{M}_{\mathrm{C}} & \mathbf{0} \\
            \mathbf{0} & \mathbf{0} \\
        \end{bmatrix},~~\mathbf{K} = \begin{bmatrix}
            \mathbf{K}_{\mathrm{C}} & \mathbf{0} \\
            \mathbf{0} & \mathbf{K}_{\mathrm{I}} \\
        \end{bmatrix},~~\mathbf{B}=\begin{bmatrix}
            \mathbf{R}_{\mathrm{C}} & -\mathbf{R}_{\mathrm{I}} \\
        \end{bmatrix},~~\mathbf{j} = \begin{bmatrix}
            \mathbf{j}_{\mathrm{C}} \\
            \mathbf{j}_{\mathrm{I}} \\
        \end{bmatrix}.
    \end{equation*}
    Finally, we choose the implicit Euler method as our time-stepping scheme, for which the iteration scheme, given a constant step size $\Delta t$, is expressed as
    \begin{equation}
        \begin{bmatrix}
            \mathbf{W} & \Delta t\mathbf{B}^{\top} \\
            \Delta t\mathbf{B} & \mathbf{0} \\
        \end{bmatrix}\begin{bmatrix}
            \mathbf{a}^{(\ell+1)} \\
            \mathbf{m}^{(\ell+1)} \\
        \end{bmatrix} = \begin{bmatrix}
            \mathbf{f}^{(\ell+1)} \\
            \mathbf{0}
        \end{bmatrix}\label{eq:euler}
    \end{equation}
    using $\mathbf{W}=\mathbf{M} + \Delta t\mathbf{K}$ and $\mathbf{f}^{(\ell+1)}=\mathbf{M}\mathbf{a}^{(\ell)} + \Delta t\mathbf{j}^{(\ell+1)}$. The superscript $\bullet^{(\ell)}$ refers to time step $t_{\ell}$, for which we employ $\ell\in\{0,\ldots,n_{\mathrm{t}}\}$ such that $t_0=0$ and $t_{n_{\mathrm{t}}}=T$.

    \section{Domain Decomposition and Tree-Cotree Gauging}\label{sec:DD+TCG}
    We employ the IETI-DP principles, as in \cite{Mally_2025ab} for magnetostatics, in \eqref{eq:euler}, which implies a splitting of $\mathbf{a}^{(\ell+1)}$ into three different parts. First, a part $\mathbf{a}^{(\ell+1)}_{\mathrm{e}}$ which is eliminated from the system (gauge or Dirichlet conditions). Then, we have the primal DOFs $\mathbf{a}^{(\ell+1)}_{\mathrm{p}}$, which are coupled strongly. The remaining DOFs are denoted as $\mathbf{a}^{(\ell+1)}_{\mathrm{r}}$. We assume that no coupling constraints between $\mathbf{a}^{(\ell+1)}_{\mathrm{p}}$ and $\mathbf{a}^{(\ell+1)}_{\mathrm{r}}$ exist. Then, we can split $\mathbf{B}$ into the diagonal blocks $\mathbf{B}_{\mathrm{rr}}$ and $\mathbf{B}_{\mathrm{pp}}$ with corresponding multipliers $\mathbf{m}_{\mathrm{r}}$ and $\mathbf{m}_{\mathrm{p}}$. The primal constraints associated with $\mathbf{m}_{\mathrm{p}}$ are eliminated by using a matrix $\mathbf{N}$ that represents the kernel of $\mathbf{B}_{\mathrm{pp}}$, i.e., $\mathbf{N}$ prescribes one value $p_k^{(\ell+1)}$ for each group of coupled primal DOFs. Consequently, we can express the strong coupling as $\mathbf{a}^{(\ell+1)}_{\mathrm{p}}=\mathbf{N}\mathbf{p}^{(\ell+1)}$ and reformulate \eqref{eq:euler} to obtain
    \begin{equation}
        \begin{bmatrix}
            \mathbf{W}_{\mathrm{rr}} & \mathbf{W}_{\mathrm{rp}}\mathbf{N} & \Delta t\mathbf{B}_{\mathrm{rr}}^{\top} \\
            \mathbf{N}^{\top}\mathbf{W}_{\mathrm{pr}} & \mathbf{N}^{\top}\mathbf{W}_{\mathrm{pp}}\mathbf{N} & \mathbf{0} \\
            \Delta t\mathbf{B}_{\mathrm{rr}} & \mathbf{0} & \mathbf{0} \\
        \end{bmatrix}\begin{bmatrix}
            \mathbf{a}^{(\ell+1)}_{\mathrm{r}} \\
            \mathbf{p}^{(\ell+1)} \\
            \mathbf{m}^{(\ell+1)}_{\mathrm{r}} \\
        \end{bmatrix} = \begin{bmatrix}
            \mathbf{f}^{(\ell+1)}_{\mathrm{r}} - \mathbf{W}_{\mathrm{re}}\mathbf{a}^{(\ell+1)}_{\mathrm{e}} \\
            \mathbf{N}^{\top}\left(\mathbf{f}^{(\ell+1)}_{\mathrm{p}} - \mathbf{W}_{\mathrm{pe}}\mathbf{a}^{(\ell+1)}_{\mathrm{e}}\right) \\
            \mathbf{0} \\
        \end{bmatrix}.
    \end{equation}
    For parallelization, the typical dual-primal approach employs two sequential Schur complements. For details, we refer to \cite{Farhat_2001aa}. Here, we focus on the DOF splitting for which we carry out a tree-cotree decomposition as in \cite{Mally_2025ab}. The basic idea is to construct a tree on the underlying (control) mesh and to eliminate the DOFs belonging to the tree to gauge the system \cite{Albanese_1988aa}. As in \cite{Mally_2025ab}, the tree is constructed globally on the wirebasket first, before it is extended into the subdomain faces and at last into the subdomain interiors. A challenge arises because only the subproblem in $\Omega_{\mathrm{I}}$ can be gauged as $\mathbf{W}_{\mathrm{I}}=\Delta t\mathbf{K}_{\mathrm{I}}$ is singular, while $\mathbf{W}_{\mathrm{C}}=\mathbf{M}_{\mathrm{C}} + \Delta t\mathbf{K}_{\mathrm{C}}$ is non-singular. In other words, we need to eliminate all tree DOFs from $\Omega_{\mathrm{I}}$, but not from $\Omega_{\mathrm{C}}$. Consequently, we need to be careful on $\Gamma$ because, to remain consistent, coupled tree DOFs from $\Omega_{\mathrm{I}}$ cannot be determined arbitrarily, but they have to take the corresponding value from $\Omega_{\mathrm{C}}$. A remedy to this is to select the tree DOFs on $\Gamma$ as primal, i.e., to select them as $\mathbf{a}^{(\ell+1)}_{\mathrm{p}}$. A sketched visualization for the tree construction and the consistent gauging is provided in Fig.~\ref{fig:tree-cotree}. Accordingly, all Dirichlet DOFs and all tree DOFs in the interior of $\Omega_{\mathrm{I}}$ are selected as $\mathbf{a}^{(\ell+1)}_{\mathrm{e}}$. Note that this elimination can entail the aforementioned homogenization procedure for inhomogeneous Dirichlet boundary conditions. All DOFs that remain (cotree on $\Gamma$ and $\Omega_{\mathrm{I}}$, all DOFs in the interior of $\Omega_{\mathrm{C}}$) are selected as $\mathbf{a}^{(\ell+1)}_{\mathrm{r}}$, which yields a non-singular $\mathbf{W}_{\mathrm{rr}}$.

    \section{Numerical Experiments}\label{sec:numers}
\begin{figure}
        \centering
        \includegraphics[width=\linewidth, trim=4.5cm 5.3cm 5cm 5.6cm, clip]{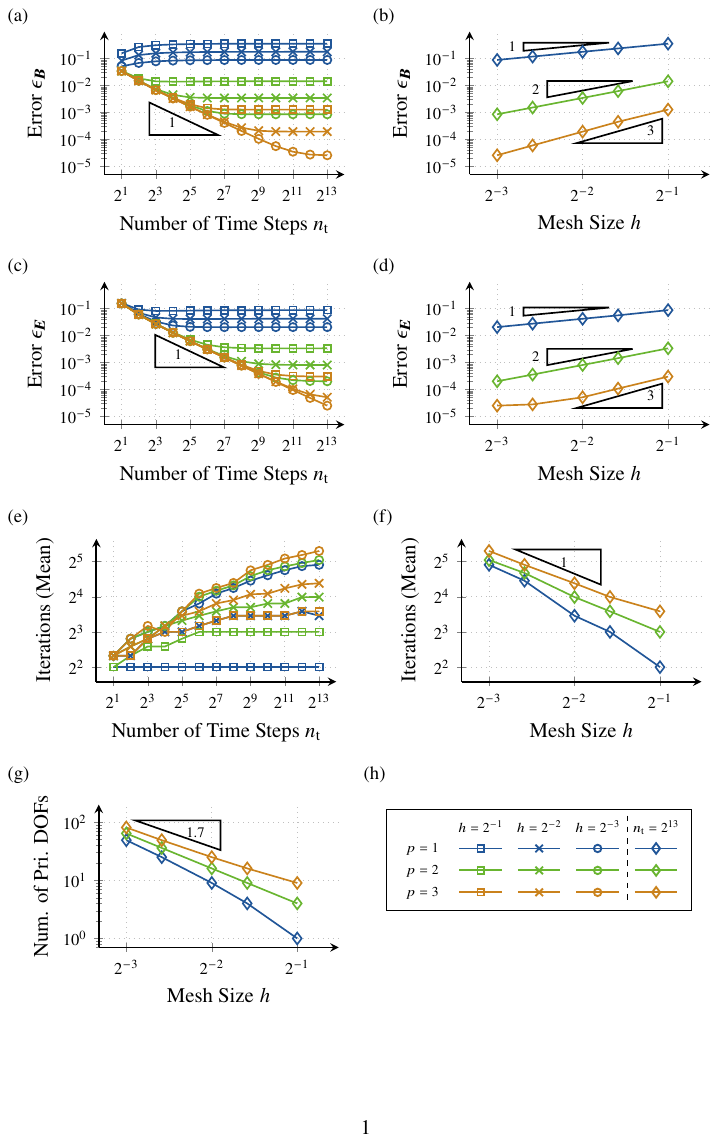}
        \caption{Measurements for numerical experiments with PCG tolerance $\epsilon=10^{-6}$.}
        \label{fig:astar_cube2}
    \end{figure}

    To construct an analytically solvable test problem, for which our implementations are available at \cite{Mally_2025ah} and based on \texttt{GeoPDEs} \cite{Vazquez_2016aa}, we prescribe the solution
    \begin{equation}
        \Afield=e^{-t}\begin{bmatrix}
            \sin(x)\cos(y)\cos(z) \\
            -2\cos(x)\sin(y)\cos(z) \\
            \cos(x)\cos(y)\sin(z) \\
        \end{bmatrix}~~\Rightarrow~~\Bfield=\Curl\Afield,~~\Efield_{\mathrm{C}} = -\ddt\Afield_{\mathrm{C}},\label{eq:problem}
    \end{equation}
    where we only use $\Efield_{\mathrm{C}}$ because it is only uniquely defined in $\Omega_{\mathrm{C}}$. From \eqref{eq:problem}, we derive appropriate homogenized boundary, source and initial terms on $\Omega=(0,1)^3$ for $t\in(0,1)$. The interface $\Gamma=\{x=0.5,~y,z\in(0,1)\}$ splits $\Omega$ into $\Omega_{\mathrm{C}}$ (for $x<0.5$) and $\Omega_{\mathrm{I}}$ (for $x>0.5$) with $\reluctivity_{\mathrm{C}}=\reluctivity_{\mathrm{I}}=1$ and $\conductivity=1$.
    Using the errors $\epsilon_{\Efield}$ for $\Efield_{\mathrm{C}}$ and $\epsilon_{\Bfield}$ for $\Bfield$ as defined in \cite[p.~141]{Acevedo_2013aa}, we expect to observe the bounds $\epsilon_{\Efield}=\mathcal{O}(h^p + n_{\mathrm{t}}^{-1})$ and $\epsilon_{\Bfield}=\mathcal{O}(h^p + n_{\mathrm{t}}^{-1})$ with mesh size $h$ and number of time steps $n_{\mathrm{t}}$. This is verified by the results in Fig.~\ref{fig:astar_cube2}a-\ref{fig:astar_cube2}d. TI methods have two main indicators for scalability. First, the number of iterations required to solve the interface problem for which we observe a dependency on both $h$ and $n_{\mathrm{t}}$ in Fig.~\ref{fig:astar_cube2}e and Fig.~\ref{fig:astar_cube2}f. Note that we employed a Dirichlet preconditioner without scaling and that we computed the mean number of iterations over all time steps. Here, we can only state that the linear increase with decreasing $h$ is not optimal for TI methods because optimal growth would be $\mathcal{O}((1 - \log(h))^2)$ as stated in \cite{Farhat_2001aa}. The second indicator, plotted in Fig.~\ref{fig:astar_cube2}g, is the number of primal DOFs which clearly depends on $h$. This behavior is not optimal for scalability but is expected because all tree DOFs on the interface are selected as primal. Therefore, the number of primal DOFs is related to the configuration of the interface mesh (of vertices and edges), which is, in turn, linked to the mesh size $h$. Our \texttt{Matlab} implementation is not parallelized, but we measured the computational time of our algorithm to give an intuition of the performance. For $n_t = 2^5$, $p=3$ and $h = 2^{-3}$, which corresponds to 7260 DOFs in $\mathbf{a}^{(\ell+1)}$, the total computational time is about $10\,\mathrm{s}$ on a laptop. For the provided example, we expect no relevant, parallel speed-up because we only consider two subdomains and the scalability indicators are not behaving optimally.

    \section{Conclusion and Outlook}
    The goal of this paper is to present a way to introduce a TI method for 3D eddy current problems in time-domain. This entails coupling of conducting $\Omega_{\mathrm{C}}$ and non-conducting $\Omega_{\mathrm{I}}$ regions consistently while enabling parallel solving of local problems in each region. To obtain this using the tree-cotree decomposition, we analyzed that selecting the tree DOFs on $\Gamma$ as primal is essential. The approach satisfies classical convergence bounds but is not optimal in regard to the typical scalability indicators without techniques or modifications that further improve performance. For future improvements of this approach, we recommend investigating scaling and deflation in the Dirichlet preconditioner or other preconditioning techniques. Furthermore, the tree DOFs on the interface need to be deselected as primal while preserving local non-singularity of subdomain system matrices.

    \begin{acknowledgement}
        The work is supported by the joint DFG/FWF Collaborative Research Centre CREATOR (DFG: Project-ID 492661287/TRR 361; FWF: 10.55776/F90) at TU Darmstadt, TU Graz and JKU Linz. Furthermore, this work has received financial support from the Consellería de Educación, Ciencia, Universidades e Formación Profesional - Xunta de Galicia (ED431C 2025/09 and ED431F 2025/03)
    \end{acknowledgement}

    \bibliography{bibtex}

@article{Acevedo_2013aa,
    author = "Acevedo, Ramiro and Loaiza, Gerardo",
    number = "17",
    pages = "111--145",
    title = "A Fully-Discrete Finite Element Approximation for the Eddy Currents Problem",
    volume = "9",
    year = "2013",
    language = "english",
    journal = "ing. cienc.",
}

@article{Albanese_1988aa,
    author = "Albanese, Raffaele and Rubinacci, Guglielmo",
    number = "7",
    pages = "457--462",
    title = "Integral Formulation for {3D} Eddy-Current Computation Using Edge Elements",
    volume = "135",
    year = "1988",
    month = "09",
    language = "english",
    journal = "{IEE} Proc. Sci. Meas. Tech.",
}

@book{Alonso-Rodriguez_2010aa,
    author = "Alonso Rodríguez, Ana and Valli, Alberto",
    publisher = "Springer",
    series = "Modeling, Simulation and Applications",
    title = "Eddy Current Approximation of {Maxwell} Equations",
    volume = "4",
    year = "2010",
    language = "english",
    address = "Heidelberg",
}

@article{Buffa_2020aa,
    author = "Buffa, Annalisa and Corno, Jacopo and de Falco, Carlo and Schöps, Sebastian and Vázquez, Rafael",
    funding = "dfg-igabem",
    number = "1",
    pages = "B80--B104",
    title = "Isogeometric Mortar Coupling for Electromagnetic Problems",
    volume = "42",
    year = "2020",
    month = "01",
    language = "english",
    journal = "{SIAM} J. Sci. Comput.",
}

@book{Cottrell_2009aa,
    author = "Cottrell, J. Austin and Hughes, Thomas Joseph Robert and Bazilevs, Yuri",
    publisher = "Wiley",
    isbn = "978-0470749098",
    title = "Isogeometric Analysis: Toward Integration of {CAD} and {FEA}",
    year = "2009",
}

@article{Emson_1988aa,
    author = "Emson, C. R. I. and Trowbridge, Charles William",
    number = "1",
    pages = "86--89",
    title = "Transient 3D eddy currents using modified magnetic vector potentials and magnetic scalar potentials",
    volume = "24",
    year = "1988",
    month = "01",
    journal = "{IEEE} Trans. Magn.",
}

@article{Farhat_2001aa,
    author = "Farhat, Charbel and Lesoinne, Michel and LeTallec, Patrick and Pierson, Kendall and Rixen, Daniel",
    number = "7",
    pages = "1523--1544",
    title = "{FETI-DP}: a dual–primal unified {FETI} method—part I: A faster alternative to the two-level {FETI} method",
    volume = "50",
    year = "2001",
    journal = "Int. J. Numer. Meth. Eng.",
}

@article{Mally_2025ab,
    author = "Mally, Mario and Kapidani, Bernard and Merkel, Melina and Schöps, Sebastian and Vázquez, Rafael",
    doi+an = "=openaccess",
    funding = "dfg-trr361,a02,c06,mercator,ernstludwig",
    ids = "Mally_2024ae",
    pages = "117737",
    related = "Mally\_2025ac",
    relatedtype = "data",
    title = "Tree-Cotree-Based Tearing and Interconnecting for {3D} Magnetostatics: A Dual-Primal Approach",
    volume = "437",
    year = "2025",
    journal = "Comput. Meth. Appl. Mech. Eng.",
}

@article{Mally_2025ah,
    author = "Mally, Mario",
    doi = "10.5281/zenodo.17237421",
    doi+an = "=openaccess",
    funding = "a02,dfg-trr361",
    howpublished = "Zenodo",
    title = "{IETIDP\_4\_MQS}",
    year = "2025",
    month = "10",
}

@book{Salon_2023aa,
    author = "Salon, Sheppard J. and Chari, M. V. K. and Ergene, Lale T. and Burow, David and DeBortoli, Mark",
    publisher = "John Wiley {\&} Sons, Ltd",
    title = "Eddy Currents: Theory, Modeling and Applications",
    year = "2023",
    language = "english",
}

@article{Vazquez_2016aa,
    author = "Vázquez, Rafael",
    number = "3",
    pages = "523--554",
    title = "A new design for the implementation of isogeometric analysis in {Octave} and {Matlab}: {GeoPDEs} 3.0",
    volume = "72",
    year = "2016",
    month = "08",
    journal = "Comput. Math. Appl.",
}

@article{Yao_2012ab,
    author = "Yao, Wang and Jin, Jian-Ming and Krein, Philip T.",
    number = "4",
    pages = "1078--1086",
    title = "A Highly Efficient Domain Decomposition Method Applied to {3-D} Finite-Element Analysis of Electromechanical and Electric Machine Problems",
    volume = "27",
    year = "2012",
    journal = "{IEEE} Trans. Energ. Convers.",
}

\end{document}